\input  amstex
\input amsppt.sty
\magnification1200
\vsize=23.5truecm
\hsize=16.0truecm
\NoBlackBoxes
\nologo

\def\rn{{\Bbb R}^n}
\def\rnp{{\Bbb R}^n_+}
\def\rnm{\Bbb R^n_-}

\def\crnp{\overline{\Bbb R}^n_+}
\def\crnm{\overline{\Bbb R}^n_-}

\def\comega{\overline\Omega }
\def\ang#1{\langle {#1} \rangle}
\def\rp{ \Bbb R_+}

\def\rmi{ \Bbb R_-}

\def\N{\Bbb N}
\def\R{\Bbb R}
\def\C{\Bbb C}

\def\ol{\overline}

\def\simto{\overset\sim\to\rightarrow}

\document

\document
\topmatter
\title
Fourier methods for fractional-order operators 
\endtitle
\author Gerd Grubb \endauthor
\affil
{Department of Mathematical Sciences, Copenhagen University,\\
Universitetsparken 5, DK-2100 Copenhagen, Denmark.\\
E-mail {\tt grubb\@math.ku.dk}}\endaffil
\rightheadtext{Fourier methods}

\endtopmatter

This is the written material for two lectures given at the conference at RIMS, Kyoto: "Harmonic Analysis and Nonlinear Partial Differential equations",
July 11-13, 2022.

The intention was to explain how methods using Fourier transformation
and complex analysis lead to sharp regularity results in the study of
fractional-order operators such as $(-\Delta )^a$, the fractional
Laplacian ($0<a<1$), and to give an overview of the results. As
required by the organizers, we start at a fairly
elementary level, introducing the role of function spaces and linear
operators. In the later text we  explain two important points in
detail, with an elementary argumentation: How
the exact solution spaces (the $a$-transmission spaces) come into the picture, and why a
locally defined Dirichlet boundary value is relevant.

Here is a small selection of the many contributors to the field:
  Blumenthal and Getoor \cite{BG59}, Vishik and Eskin '60s (presented in \cite{E81}), Hoh and Jacob
\cite{HJ96},  Kulczycki \cite{K97}, Chen and Song \cite{CS98},
Jakubowski \cite{J02}, Bogdan, Burdzy and Chen \cite{03}, Cont and
Tankov \cite{04}, Caffarelli and Silvestre \cite{07}, Gonzales, Mazzeo and Sire \cite{12},
Ros-Oton and Serra \cite{RS14}, Grubb \cite{G15}, Abatangelo \cite{A15}, Felsinger, Kassmann and Voigt
\cite{FKV15}, Bonforte, Sire and Vazquez \cite{BSV15}, Dipierro, Ros-Oton and Valdinoci
\cite{DRV17}, Dyda, Kuznetzov and Kvasnicki \cite{DKK17}, Abatangelo, Jarohs and Saldana \cite{AJS18}, Chan, Gomez-Castro and Vazquez \cite{CGV21},
  Borthagaray and Nochetto \cite{BN23}.
Besides these works listed at the end, we shall only list
 the papers that are directly referred to in the text. Many more
 references are given in the works.

\medskip

\noindent {\bf Plan of the lectures:}
\smallskip
 
\noindent  { 1. The homogeneous Dirichlet problem:}

  {  1.1.  Introduction, the Fourier transform.}

  {  1.2. The fractional Laplacian.}

   {  1.3. Model Dirichlet problems. (Detailed)}

  {  1.4. The Dirichlet problem for curved domains.} 
  \smallskip

\noindent  {  2. Further developments:}
  
  {  2.1. Evolution problems and resolvents.}

 {  2.2. Motivation for local nonhomogeneous boundary conditions. (Detailed)}

{2.3. Nonhomogeneous Dirichlet conditions over curved
domains.}

 {  2.4. Integration by parts,  Green's formula.}

References

\head 1. The homogeneous Dirichlet problem \endhead 

\subhead
1.1 Introduction, the Fourier transform \endsubhead

We start by 
recalling the basic notions of {\it function spaces} and {\it operators}:
\medskip

{\bf Function spaces.}
When $f(x)$ is a function on Euclidean space $\Bbb R^n$, with
points denoted $x=(x_1,x_2,\dots,x_n)$, differentiation with respect
to each variable $x_k$ gives the {\it partial derivative}
$$
\frac{\partial f(x_1,\dots,x_n)}{\partial x_k},\text{ also denoted }\partial_kf.
$$
Here are some examples of spaces of functions with derivatives:
\smallskip

$\bullet$ $C^0(\Bbb R^n)$ consists of the bounded continuous
functions on $\Bbb R^n$.

\smallskip

$\bullet$ $C^m(\Bbb R^n)$ ($m\in{\Bbb N}$) consists of those bounded continuous
functions $f$ that allow taking  partial derivatives up
to $m$ times giving bounded continuous function.

 \smallskip
 $\bullet$ $L_2(\Bbb R^n)$ consists of functions $f$
such that $\int_{\rn} |f(x)|^2 \, dx$ exists. (Here one uses
Lebesgue's measure theory, identifying functions that coincide outside
 a null-set.)

\smallskip
 
$\bullet$ {\bf Sobolev spaces} $H^m(\rn)$ consist of the functions $f$ in
$L_2(\rn )$ that have partial derivatives up to order $m$ in $L_2(\rn
)$ (in
a generalized sense).

\smallskip
 Each of these function spaces is a linear infinite dimensional vector
space (when $f$ and $g$ are there, the sum $c_1f(x)+c_2g(x)$ is
likewise there). They are normed spaces.

\medskip 
  
{\bf Operators.}
Linear operators are mappings from one function space to another,
preserving the vector space structure.

For example,
$\partial_k$ defines a  linear  operator going
from $C^m(\Bbb R^n)$ to $C^{m-1}(\Bbb R^n)$, and from
$H^m(\rn)$ to  $H^{m-1}(\rn)$, when $m\ge 1$.
More generally, a {\bf
  partial differential operator} is a
sum of composed derivatives
multiplied by functions, $A=\sum_{|\alpha |\le k}a_\alpha (x)\partial  ^\alpha $. 
(Here we use the multi-index notation: Let $\alpha =(\alpha _1,\dots,\alpha _n)\in
\Bbb N_0^n$; $\N_0=\{0,1,2,\dots\}$. Then $x^\alpha =x_1^{\alpha _1}\cdots x_n^{\alpha _n}$,
$\partial^\alpha =\partial_1^{\alpha _1}\cdots \partial_n^{\alpha
  _n}$, $|\alpha |=\alpha _1+\dots+\alpha _n$.)
  $A$ goes from  $C^m(\Bbb R^n)$ to $C^{m-k}(\Bbb R^n)$, and from
$H^m(\rn)$ to  $H^{m-k}(\rn)$, when $m\ge k$ and the coefficients
$a_\alpha $ are smooth and bounded with bounded derivatives.

A very important
example is {\bf the Laplace operator}
$$
\Delta \colon u\mapsto \Delta u=\partial_1^2u+\dots+\partial _n^2u.
$$

It enters in three basic equations (two of them has an extra variable $t$): 
$$
\aligned
-\Delta u(x)&= f(x)\text{ on }\Omega , \text{ the Laplace equation,}\\
\partial_tu(x,t)-\Delta u(x,t)&= f(x,t)\text{ on }\Omega \times
\Bbb R, \text{ the heat equation,}\\
\partial_t^2u(x,t)-\Delta u(x,t)&= f(x,t)\text{ on }\Omega \times
\Bbb R,\text{ the wave equation,}
\endaligned
$$  
describing physical problems. Here $\Omega $ is an open subset of $\rn$, and
one wants to  {\it find solutions} $u$ for given $f$.
For example, the heat equation describes how the temperature develops
in a container. And it is also used in financial theory, wrapped up in
a stochastic formulation.

Another type of examples of operators are {\bf integral operators},
such as
$$
(\Cal K u)(x)=\int _{\rn}  K(x,y)u(y)\, dy.
$$
Differential equations often have integral operators as
solution operators.

A very important special operator is the  {\bf Fourier transformation}
$\Cal F$, it is an integral operator:
$$\Cal F
u= \hat
u(\xi )=
\int_{{\Bbb R}^n}e^{-ix\cdot \xi }u(x)\, dx.
$$
  It is invertible (in fact isometric times a constant) from the space
$L_2(\Bbb R^n)$ onto $L_2(\Bbb R^n)$, and the inverse operator
looks similar: $\Cal F^{-1}v=
c\int_{{\Bbb R}^n}e^{+ix\cdot \xi }v(\xi )\, d\xi 
$, $c=(2\pi )^{-n}$.
It has been used much in physics and
mathematics, and the mathematical rigor was perfected with 
Schwartz'
Distribution Theory around 1950, defining the rapidly decreasing functions
$\Cal S(\rn)$ and temperate distributions $\Cal S'(\rn)$.  (A detailed presentation is given
e.g.\ in \cite{G09}.)
 
The success of $\Cal F$
comes from the
fact that {\it it turns differential operators into multiplication operators}: 
The differential operator $\partial_k$ is turned into multiplication
by $i\xi _k$:
$$
\Cal F(\partial_ku)=i\xi _k\hat u(\xi )\qquad \text{ (here }i=\sqrt{-1}).
$$
For example, $\Cal F(\Delta u)=-(\xi _1^2+\dots+\xi _n^2)\hat
u(\xi )=-|\xi |^2\hat u$,
and therefore the equation $
-\Delta u=f$ is turned into $ |\xi |^2\hat u=\hat f$.

For a particularly simple example, consider the operator $1-\Delta $
on $\rn$.
$$
(1-\Delta )u=f\text{ is transformed to }
(1+|\xi |^2)\hat u=\hat f,$$
which has the unique solution $u=\Cal F^{-1}(\tfrac1{1+|\xi
  |^2}\hat f).$

\medskip

  {\bf Pseudodifferential operators.}
Now we  generalize the above idea:  Take a function $p(\xi )$, the
{\it symbol},  and {\it define} the {\it pseudodifferential
operator} ($\psi $do) $P=\operatorname{Op}(p)$ by
$$
\operatorname{Op}(p)(u)=\Cal F^{-1}(p(\xi )\hat u(\xi ))=
\Cal F^{-1}p(\xi )\Cal F u.
$$
Then if $p$ has an inverse $1/p$,
$\operatorname{Op}(p)\operatorname{Op}(1/p)=\operatorname{Op}(p\cdot1/p)=I
$, so   $P=\operatorname{Op}(p)$ has the
inverse $P^{-1}=\operatorname{Op}(1/p)$. This is the simple basic idea.

We often need to let the  symbol  $p$ depend on $x$ also. This
is natural, since differential operators in general have $x$-dependent
coefficients, but it gives more difficult composition rules.
The definition is  
$$ 
\operatorname{Op}(p(x,\xi ))u(x)=\tfrac1{(2\pi )^n}\int_{\rn} e^{ix\cdot\xi
                                  }p(x,\xi )\hat u(\xi )\, d\xi
$$ 
under suitable requirements on  $p(x,\xi )$. We say
that $p$ is {\it of order} $m$ when $\partial_x^\beta \partial_\xi ^\alpha p$
is $O((1+|\xi |)^{m-|\alpha |})$ for all multi-indices $\alpha ,\beta $.

For $x$-independent symbols there is the simple composition rule
$$
\operatorname{Op}(a(\xi ))\operatorname{Op}(b(\xi
))u=\Cal F^{-1}a(\xi )\Cal F\Cal F^{-1}b(\xi )\Cal Fu=\operatorname{Op}(a(\xi )b(\xi ))u;
$$
in other words, the symbol of
  $\operatorname{Op}(a)\operatorname{Op}(b)$ is $ab$.
  For $x$-dependent symbols there is, just like for differential
  operators, a more complicated composition formula with lower-order terms.
$$
  \operatorname{Op}(a(x,\xi )) \operatorname{Op}(b(x,\xi ))=
\operatorname{Op}(a(x,\xi )b(x,\xi ))+\Cal R,
$$
where the order of $\Cal R$ is 1 step lower  than that of $ab$
($\Cal R$ can be described in more detail).

The case of $x$-independent symbols can often be used as a model for
the general case.

\smallskip
The theory was built up in the 1960's (by Kohn and Nirenberg, H\"ormander, Seeley, with preceding insights by Mihlin,
Calderon, Zygmund and others), and further developed through the rest
of the century and beyond.

\subhead 1.2 The fractional Laplacian \endsubhead

The operator we shall be concerned with here is the fractional Laplacian
$(-\Delta )^a$, $0<a<1$. 
It can be defined by spectral theory in functional analysis, since
$-\Delta $ is a selfadjoint nonnegative operator (unbounded) in the Hilbert space
$L_2(\rn)$.  It is currently of 
great interest in probability theory and
finance, and also in mathematical physics and differential geometry.

Structurally, it
is a {\bf pseudodifferential operator},
$$
(-\Delta )^au=\operatorname{Op}(|\xi |^{2a})u. \tag1.1
$$
  
It can also be written as a {\bf singular integral operator}:
$$
   (-\Delta )^au(x)=c_{n,a}PV\int_{{\Bbb
    R}^n}\frac{u(x)-u(x+y)}{|y|^{n+2a}}\,dy; \tag 1.2
$$
here $c_{n,a}|y|^{-n-2a}=\Cal F^{-1}|\xi
|^{2a}$, and PV stands for ``principal value''.

%

Formula
(1.1) has natural generalizations to $x$-dependent symbols $p(x,\xi )$,
allowing ``variable-coefficient'' operators.

Formula (1.2) is often used in probability and nonlinear analysis, 
with generalizations to expressions with other
kernel functions than $|y|^{-n-2a}$, e.g.\ $|y|^{-n-2a}K(y/|y|)$, $K$
positive, and {\it even:} $K(-y)=K(y)$ (possibly with less smoothness). They
generate L\'evy processes. Here calculations are often made
considering only real functions, whereas the Fourier transform of course
involves complex functions.

In contrast to $-\Delta $, $(-\Delta )^a$ is a {\it nonlocal}
operator on $\rn$: When $u=0$ in an open set $\omega $, then $\Delta u=0$ on
$\omega $
 but usually $(-\Delta )^au \ne 0$ there;
this gives substantial difficulties.
 To study functions $u$  on a given open subset
  $\Omega $ of $\rn$, we can define
  $(-\Delta )^au$ by letting $u$ be zero on
  $\rn\setminus \comega $ (i.e., $\operatorname{supp}u\subset\overline\Omega
  $), and map it to $r^+(-\Delta )^au$, where $r^+$ denotes restriction to $\Omega $.

The {\it homogeneous Dirichlet problem} for $P=(-\Delta )^a$
on $\Omega $ is then defined
as follows: For a given function $f$ on $\Omega $,  find a function
$u$ on $\rn$ such that
$$
r^+Pu= f \text{ in }\Omega ,\quad \operatorname{supp}u\subset
\overline\Omega .\tag1.3
$$
 
We now need to introduce Sobolev spaces over $\Omega $. Denote $\ang\xi
=(1+|\xi |^2)^{\frac12}$; then 
$$
\aligned
H^s(\rn)&=\{u\in \Cal S'(\Bbb R^n)\mid \langle{\xi }\rangle^s\hat
u\in L _2(\Bbb R^n)\},\text{ with norm }\|\Cal F^{-1}(\ang{\xi }^s\hat u)\|_{L_2},\\
  \overline H ^s(\Omega  )&=r^+H^s (\Bbb R^n),\text{ the {\it
  restricted} space,}\\
\dot H^s (\comega )&=\{u\in H^s (\Bbb R^n)\mid
\operatorname{supp}u\subset \comega  \},\text{ the {\it supported} space},
\endaligned
$$
for $s\in\Bbb R$. 
(The dot and overline notation
stems from H\"ormander's books.)
When $\Omega $ is suitably regular, $\overline H^s(\Omega )$ and $\dot H^{-s}(\comega)$ are dual
spaces, with a duality consistent with the $L_2$-scalar
product. The space $C_0^\infty (\Omega )$ of smooth functions with
compact support in $\Omega $ is dense in $\dot H^s(\comega)$ for all $s$. The space $\overline H^s(\Omega )$ coincides with
$\dot H^{s}(\comega)$ when $|s|<\frac12$.

\medskip

There is a variational formulation that gives unique solvability of
(1.3) in low-order Sobolev spaces. Let $P=(-\Delta )^a$, $0<a<1$. First note
that $P$ maps $ H^{s}(\rn)\to H^{s-2a}(\rn)$ continuously, since 
$$
\|Pu\|^2_{H^{s-2a}}=c\int(\ang \xi ^{s-2a}|\xi |^{2a}|\hat u|)^2d\xi \le c\int(\ang\xi ^{s}|\hat u|)^2d\xi =\|u\|^2_{H^s}.
$$
In particular, $P\colon H^{a}(\rn)\to H^{-a}(\rn)$, and hence $r^+ P\colon \dot H^{a}(\comega)\to
\overline H^{-a}(\Omega )$.

The sesquilinear form $Q_0$ on $\dot H^a(\comega )$ obtained by closure of
$$
Q_0(u,v)=\int_{\Omega
}Pu\,\bar v\,dx\text{ for }u,v\in C_0^\infty (\Omega ),
$$
satisfies $Q_0(u,u)\ge
0$, and equals $\ang {r^+Pu,v}_{\overline H^{\,-a},\dot H^a}$.
 The $L_2$ Dirichlet realization
$P_D$ is defined as the operator acting like $r^+P$ and having domain 
$$
D(P_D)=\{u\in \dot H^a(\comega)\mid r^+Pu\in L_2(\Omega )\};
$$
this operator is selfadjoint in $L_2(\Omega )$ and $\ge 0$.
 
When $\Omega $ is bounded, 
there is a Poincar\'e
inequality assuring that $P_D$ has positive lower
bound, hence it is {\it bijective} from $D(P_D)$ onto $L_2(\Omega )$.

So there {\it exists} a bijective solution operator for (1.3), but what
more can we say about $u$, when $f\in L_2(\Omega )$, or lies in better spaces $\overline
H^s(\Omega )$?

  It has been known since the 1960's that
  $D(P_D)=\dot H^{2a}(\comega)$ if $a<\frac12$, and
  $D(P_D)\subset \dot H^{a+\frac12-\varepsilon }(\comega)$ 
  if $a\ge\frac12$, \cite{E81}, but more precise information has
 been obtained in recent years.

 The new knowledge is that under some regularity assumptions, the solution $u$ has a
 factor $d^a$,  where
   $d(x)=\operatorname{dist}(x,\partial\Omega )$. We shall quote the
  detailed results later, but will now show  how the factor $d^a$ comes
   in via Fourier transformation methods.
For this, we recall some important formulas for the Fourier transform from functions of
$x_n\in{\Bbb R}$ to functions of $\xi _n\in{\Bbb R}$:
Denoting $1|_{\rp}=H(x_n)$ (the Heaviside function), we have for $\sigma >0$, $a> -1$,
$$
\aligned
{\Cal F}_{x_n\to \xi _n }(H(x_n)e^{-\sigma x_n})&=\frac{1}{\sigma
  +i\xi _n },\\
    {\Cal F}_{x_n\to \xi _n }(H(x_n)x_n^ae^{-\sigma x_n})&=\frac{c}{(\sigma
  +i\xi _n )^{a+1}},\quad c=\Gamma
  (a+1).
  \endaligned \tag1.4
$$  
The complex number $\sigma +i\xi _n $ has real part $\sigma >0$, so its
noninteger powers (defined with a cut along the negative axis $\rmi$)
make good sense. The first formula is elementary; proofs of the second
formula are found e.g.\ in Schwartz \cite{S61, (V,1;44)} and in the lines
after Example 7.1.17 in \cite{H83} (with different conventions).
 (Using (1.4), we can avoid going in detail with homogeneous
 distributions, limits for $\sigma  \to 0$.)

     \subhead 1.3 Model Dirichlet problems\endsubhead

We shall now study the Dirichlet problem (1.3) in the simplest
  possible case where $P$ is the invertible $\psi $do $(1-\Delta )^a$
  and $\Omega =\rnp$.
  \smallskip
  
  {\bf Example 1.}   First we make some remarks on the Dirichlet problem for $1-\Delta $ on
     $\rnp$ ($=\{(x',x_n)\mid x_n>0\}$;
     $x'=(x_1,\dots,x_{n-1}\}$), denoting $u(x',0)=\gamma _0u$:
$$
(1-\Delta )u=0\text{ on }\rnp,\quad \gamma _0u=\varphi \text{ on
}\Bbb R^{n-1}. \tag1.5
$$  
Fourier transformation in $x '$ turns the operator into $1+|\xi
'|^2-\partial_n^2$, so (1.5) becomes an ODE problem for each $\xi '$:
$$
(\ang{\xi '}^2-\partial_n^2)\acute u(\xi ',x_n)=0\text{ on }\rnp,
\quad \acute u(\xi ',0)=\hat\varphi (\xi ').
$$  
This has  the unique bounded solution $\acute u(\xi
',x_n)=\hat\varphi (\xi ')e^{-\ang{\xi '}x_n}$ on  $\rp$.  
Inverse
Fourier transformation from $\xi '$ to $x'$ gives that (1.5) is solved by $u=K_0\varphi $, where
$K_0$ is the {\it Poisson operator} defined by
$$
K_0\varphi 
=\Cal F^{-1}_{\xi '\to x'}\bigl(H(x_n)e^{-\ang{\xi '}x_n}\hat \varphi (\xi ')\bigr)=\Cal F^{-1}_{\xi \to x}\bigl(\tfrac1{\ang{\xi '}+i\xi _n}\hat \varphi
(\xi ')\bigr),
$$
using (1.4).  
Here $K_0$ maps continuously $
K_0\colon H^{s-\frac12}(\Bbb R^{n-1})\to e^+ \overline H^s(\rnp)$, for
all $s\in \R$.\footnote{Strictly speaking, the standard Poisson
operator is $r^+K_0$, but this distinction is often left out, since
$K_0\varphi $ is $0$ on $\rnm$. More
comments in \cite{G14, (A-13)-(A-14)}, \cite{G19, Remark 3.2}.}

\medskip

{\bf Example 2.} Now turn to the model Dirichlet problem for the fractional Laplacian
(recall $0<a<1$):
  $$
r^+(1-\Delta )^au=f\text{ on }\rnp,\quad \operatorname{supp}u\subset
\crnp.\tag 1.6
$$  
The variational solution method applies straightforwardly to
$P=(1-\Delta )^a$, showing that for $f\in L_2(\rnp)$ there is a unique
solution $u\in \dot H^a(\crnp)$. It will now be examined.

The symbol $p(\xi )=(1+|\xi |^2)^a$ of $P$ has the factorization: 
$$
(1+|\xi |^2)^{a}=(\ang{\xi
  '}^2+\xi _n^2)^a=(\ang{\xi '}-i\xi _n)^a(\ang{\xi '}+i\xi _n)^a.
$$
 
Introduce  for general $t\in\Bbb R$ the {\it order-reducing operators}: 
$$
\Xi
_\pm^t=\operatorname{Op}((\ang{\xi '}\pm i\xi _n)^t).
$$
 They are  invertible,   mapping for
all $s\in \Bbb R$:
$$
\Xi _\pm^t\colon  H^{s}(\rn)\simto H^{s-t}(\rn),\text{ with
  inverse }\Xi _\pm^{-t}.
$$  
Using these, $ P$ has the factorization
$$
(1-\Delta )^a=\Xi _-^a\,\Xi _+^a,\text{ with inverse }(1-\Delta
)^{-a}=\Xi _+^{-a}\,\Xi _-^{-a}.
$$

 The operators $\Xi _\pm^t$ have special roles relative to $\rnp$. The {\it plus-family}  $\Xi
 _+^t$ has
 symbols that extend analytically in $\xi _n$ to the lower complex
 halfplane $\Bbb C_-=\{\operatorname{Im}\xi _n<0\}$; then by the
 Paley-Wiener theorem, $\Xi ^t_+$ {\it preserves support in
   $\crnp$}. The inverse is $\Xi _+^{-t}$. Thus  for all $s\in\Bbb
 R$,
 $$\Xi _+^t\colon \dot H^{s}(\crnp)\simto \dot
 H^{s-t}(\crnp). \tag1.7
 $$  
 The {\it minus-family} $\Xi _-^t$ behaves in a similar way with respect to
 $\crnm$. Since $\Xi _-^t=(\Xi _+^t)^*$, we have moreover, in view of the
 duality between  $\dot H^{s}(\crnp)$ and $\overline H^{\,-s}(\rnp)$, that
 $$
 r^+\Xi _-^te^+\colon \overline H^{s}(\rnp)\simto \overline
H^{s-t}(\rnp), \tag1.8
$$
with inverse $(r^+\Xi _-^te^+)^{-1}=r^+\Xi _-^{-t}e^+$; here $e^+$
indicates ``extension by zero'' from $\rnp$ to $\rn$. (For negative
values of $s$, there is a distributional interpretation of (1.8).)

Note also that for $v\in L_2(\rn)$, $v=e^+r^+v+e^-r^-v$, by the identification of 
  $L_2(\rn)$ with $e^+L_2(\rnp)\,\dot + \,e^-L_2(\rnm)$. 

\medskip

Let $u\in \dot
  H^a(\crnp)$; then $\Xi _+^a$ maps it into $e^+L_2(\rnp)$.
Now since $P=\Xi _-^a\Xi _+^a$, we may write
$$
 r^+Pu=r^+\Xi _-^a\Xi _+^au=r^+\Xi _-^a(e^+r^++e^-r^-)\Xi _+^au=r^+\Xi _-^ae^+r^+\Xi _+^au,
$$
where we used that $r^-\Xi _+^au=0$. In a diagram, $r^+P$ is the
composition
$$
\dot H^{a}(\crnp)\overset{r^+\Xi  _+^{a}} \to \longrightarrow
L_2(\rnp)\overset{r^+\Xi _-^{a}e^+} \to\longrightarrow \overline
H^{\,-a}(\rnp).
$$  
Here both factors are bijections, in view of (1.7) and (1.8).  
Hence the inverse $R$, the solution operator for (1.6) with $f\in
\overline H^{\,-a}(\rnp)$, is the composed operator $R=\Xi _+^{-a}e^+(r^+\Xi _-^{-a}e^+)$;
$$
\overline H^{\,-a}(\rnp)\overset {r^+\Xi _-^{-a}e^+}\to \longrightarrow L_2(\rnp)\overset{\Xi
  _+^{-a}e^+}\to \longrightarrow \dot H^a(\crnp). \tag1.9
$$

  To find the solution of (1.6) with $f\in L_2(\rnp)$, we restrict the
  operator $R$ to $ L_2(\rnp)$; this is expressed in the diagram
  $$
  L_2(\rnp)\overset {r^+\Xi _-^{-a}e^+}\to \longrightarrow \overline H^a(\rnp)\overset{\Xi
  _+^{-a}e^+}\to \longrightarrow \Xi_+ ^{-a}e^+\overline
  H^a(\rnp)=D(P_D). \tag1.10
$$  
Property (1.8) is used in the first mapping, but in the second mapping
there can be a mismatch; property (1.7) may not be used. The space at the
right end is the so-called {\it $a$-transmission space}
$H^{a(2a)}(\crnp)$, generally defined for $t>a-\frac12$ by
$$
H^{a(t)}(\crnp)= \Xi _+^{-a}e^+\overline
H^{t-a}(\rnp).
\tag1.11
$$

The idea can also be applied starting with $f$ given in a space $\overline H^s(\rnp)$
with $s\ge -a$:
 $$
 \overline H^s(\rnp)\overset {r^+\Xi _-^{-a}e^+}\to \longrightarrow \overline H^{s+a}(\rnp)\overset{\Xi
  _+^{-a}e^+} \to\longrightarrow \Xi_+ ^{-a}e^+\overline
  H^{s+a}(\rnp)\equiv H^{a(s+2a)}(\crnp), \tag1.12
$$
with bijective mappings. This proves  
 \proclaim{Theorem 1.1}
 Let $s\ge -a$. The solution $u$ of the
  Dirichlet problem {\rm (1.6)} satisfies
$$
f\in \overline H^s(\rnp)\iff u\in H^{a(s+2a)}(\crnp).
$$
In particular, $D(P_D)=H^{a(2a)}(\crnp)$; the case $s=0$.
\endproclaim

What are these transmission spaces? Note that they decrease with
 increasing $s$, and $ H^{a(a)}(\crnp)= \dot H^{a}(\crnp)$.
 We will study $D(P_D)= \Xi_+ ^{-a}e^+\overline H^a(\rnp)$ $\equiv
 H^{a(2a)}(\crnp)$ more closely:   

For $a<\frac12$, $\overline H^a(\rnp)$ identifies with $\dot
H^a(\crnp)$ so here (1.7) can be applied and gives that $D(P_D)=\dot
H^{2a}(\crnp)$.

For $\frac12\le a<1$, we are in a new situation. Since $\overline
H^a(\rnp)\subset \ol H^{\frac12-\varepsilon }(\rnp)=\dot H^{\frac12-\varepsilon }(\crnp)$
(small
$\varepsilon >0$), we have at least $D(P_D)\subset  \dot
H^{a+\frac12-\varepsilon }(\crnp)$, by (1.7). For $a=\frac12$, the
sharpest information is $H^{\frac12(1)}(\crnp)= \Xi _+^{-\frac12}e^+\overline
H^{\frac12}(\rnp)$. For $a>\frac12$ we can analyze more:

  Formula (1.4) with $\sigma =\ang{\xi '}$ shows:
$$
{\Cal F}^{-1}_{\xi _n\to x_n}\bigl(\tfrac1{(\ang{\xi '} +i\xi _n)^{a+1}}\bigr)=\tfrac
1{\Gamma (a+1)}H(x_n)x_n^ae^{-\ang{\xi '} x_n}.\tag 1.13
$$  
Recall from Example 1 the Poisson operator $K_0\colon \varphi \mapsto u=\Cal F^{-1}_{\xi \to x}\bigl(\tfrac1{\ang{\xi '}+i\xi _n}\hat \varphi
(\xi ')\bigr)$ solving
$(1-\Delta )u=0$ on $\rnp$, $\gamma _0u=\varphi $. We can show:

\proclaim{Theorem 1.2}  $1^\circ$ For $\frac12<a<1$, $u\in H^{a(2a)}(\crnp)$
  if and only if 
$$
u=v+w,\text{ where }v\in \dot H^{2a}(\crnp), \; w=
x_n^aK_0\psi \text{ with } \psi \in H^{a-\frac12}({\Bbb R}^{n-1}).\tag 1.14
$$ 

$2^\circ$ With $\Cal S(\crnp)=r^+\Cal  S(\rn)$, $x_n^a\Cal S(\crnp)$ is dense in
$H^{a(t)}(\crnp)$. The mapping  $\gamma
_0^a\colon u\mapsto \gamma _0(u/x_n^a)$ from $x_n^a\Cal S(\crnp)$ to
$\Cal S(\R^{n-1})$ extends by continuity 
  to a mapping (when $t>a+\tfrac12$),
$$
\gamma _0^a\colon H^{a(t)}(\crnp)\to H^{t-a-\frac12}(\Bbb
R^{n-1}).
$$  
Here $\psi $ in  {\rm (1.14)} equals $\gamma ^a _0u$. 
\endproclaim
 
\demo{Proof}
$1^\circ$. Let $u\in H^{a(2a)}(\crnp)$, that is, $u=\Xi _+^{-a}e^+f$
for some $f\in \overline H^{a}(\rnp)$.
    Since $a>\frac12$, $f$ has a boundary value $\varphi
    =\gamma _0f\in H^{a-\frac12}(\R^{n-1})$, and there is a  decomposition $f=g+h$,
    where $$
    h=K_0\varphi \in
   K_0H^{a-\frac12}({\Bbb R}^{n-1}),\quad  g =f-K_0\varphi\in
 \dot H^{a}(\crnp),
$$
since $\gamma _0g=\gamma _0f-\varphi =0$.
Going back to $u$ by applying $\Xi _+^{-a}$, we find 
$$
u=v+w,\quad v=\Xi _+^{-a}g, \quad w=\Xi _+^{-a}K_0\varphi .
$$  
By the mapping property (1.7), $v\in \dot H^{2a}(\crnp)$. For $w$, we find:
$$
\aligned
w&=\Xi _+^{-a}K_0\varphi=\Cal F^{-1}_{\xi \to x}\bigl(\tfrac1{(\ang{\xi '}+i\xi
                       _n)^a} \tfrac1{\ang{\xi '}+i\xi _n}\hat\varphi (\xi ')\bigr)\\
&=\Cal F^{-1}_{\xi \to x}\bigl(\tfrac1{(\ang{\xi '}+i\xi
  _n)^{a+1}} \hat\varphi (\xi ')\bigr)\\
 &= \tfrac 1{\Gamma (a+1)}\Cal F^{-1}_{\xi '\to x'}\bigl(x_n^ae^{-\ang{\xi '}x_n}H(x_n) \hat\varphi (\xi ')\bigr)  = \tfrac 1{\Gamma (a+1)}x_n^aK_0\varphi ,
\endaligned
$$ 
where we have used formula (1.13).  
This shows  the asserted form of $w$, with $\psi =\tfrac 1{\Gamma
(a+1)}\varphi $. (We have omitted a constant entering in the
definition of $\gamma _0^a$ in the literature.)

$2^\circ$. The properties of $\gamma _0^a$  are known from \cite{G15}
(with the notation $\gamma _{0,a}u=$ \linebreak$\Gamma (a+1)\gamma _0(u/x_n^a)$); we shall just indicate a quick way
to obtain the mentioned statements. The space $\Cal
  S(\crnp)=r^+\Cal
  S(\rn)$  is dense in $\ol H^r(\rnp)$ for all $r\in\R$. There is an
  elementary proof in \cite{G21, Sect.\ 6}  of the identity
$$
e^+x_n^a\Cal S(\crnp)=\Xi _+^{-a}  e^+\Cal S(\crnp),\tag 1.15
$$
(based on  Taylor expansion in $x_n$ of
$e^{-\ang{\xi '}x_n}\Cal F_{x'\to\xi '}u$ and formulas like (1.13)); here when $u\in e^+x_n^a\Cal S(\crnp)$,
$$
\gamma _0(u/x_n^a)=c^{-1}\gamma _0(\Xi _+^au), \quad c=\Gamma (a+1).\tag 1.16
$$
(Lemma 6.1 and (6.9) in \cite{G21}.)

By (1.15), the denseness of $\Cal
S(\crnp)$ in $\ol H^{t-a}(\rnp)$ implies the denseness of  $ e^+x_n^a\Cal
S(\crnp)$ in $\Xi _+^{-a}e^+\ol H^{t-a}(\rnp)=H^{a(t)}(\crnp)$. It is
well-known that for
$t-a>\frac12$, $\gamma _0$ defined on $\Cal S(\crnp)$ extends by
continuity to the  map $\gamma _0\colon \ol H^{t-a}(\rnp)\to
H^{t-a-\frac12}(\R^{n-1})$. Then the map $u\to \gamma _0(u/x_n^a)$
from $x_n^ae^+\Cal S(\crnp)$ to $\Cal S(\R^{n-1})$
extends by continuity to a map from $H^{a(t)}(\crnp)$ to
$H^{t-a-\frac12}(\R^{n-1})$ in view of (1.16). \qed
\enddemo

 Summing up, we conclude that for $\frac12<a<1$, $D(P_D)=H^{a(2a)}(\crnp)$
  is the set of functions $u$ of the form
$$
u=v+x_n^aK_0\psi,\text{ where }\psi =\gamma _0(u/x_n^a);\tag 1.17
$$
here $v$ and $\psi $ run through $ \dot H^{2a}(\crnp)$ resp.\
$H^{a-\frac12}({\Bbb R}^{n-1})$.

  \subhead 1.4 The Dirichlet problem for curved domains \endsubhead

 For general domains $\Omega $, we shall list some recent regularity and solvability results
  in a brief formulation. First we recall the definitions of
  some more function spaces:
\smallskip

$\bullet$ The Bessel-potential spaces $H^s_q(\rn)$, $1<q<\infty $,
$s\in{\Bbb R}$, extend the Sobolev spaces $H^s(\rn)$ to $q\ne 2$:
$$
H^s_q(\Bbb R^n)
=\{u\in \Cal S'(\Bbb R^n)\mid \Cal F^{-1}(\langle{\xi }\rangle^s\hat
u)\in L_q(\Bbb R^n)\}.
$$

\smallskip

$\bullet $ The H\"older-Zygmund spaces $C^s_*(\rn)$, $s\in\Bbb R$,
generalize the H\"older spaces $C^s(\rn)$ with $s\in
\rp\setminus\Bbb N$,
to all $s$. (For $0<\sigma <1$, $u\in C^{k+\sigma }(\rn)$ with $k\in
\N_0$ and $0<\sigma <1$ when $u$ and
its derivatives up to order $k$ satisfy
$|u(x)|+\frac{|u(x)-u(y)|}{|x-y|^\sigma }\le C$ on $\rn$.)
The cases $s\in \Bbb N$ are interpolation spaces between
noninteger cases.

For an open subset $\Omega \subset \rn$, we define the scales of 
{\it restricted} spaces:
$$
\overline H_q^s(\Omega )=r^+H^s_q(\Bbb R^n), \quad  \overline C_*^s(\Omega )=r^+C_*^s(\Bbb R^n),
$$
and the scales of {\it supported} spaces:
$$
\dot H^s_q(\overline\Omega )=\{u\in H^s_q(\Bbb R^n)\mid
\operatorname{supp}u\subset \overline\Omega  \}, \quad \dot C^s_*(\overline\Omega )=\{u\in C_*^s(\Bbb R^n)\mid
\operatorname{supp}u\subset \overline\Omega  \}.
$$
 
For $\Omega =\rnp$, we define the {\bf $a$-transmission spaces}
$H_q^{a(t)}(\crnp)$ and  $C^{a(t)}_*(\crnp)$ as follows ($\frac1{q'}=1-\frac1q$): $$
\aligned
H_q^{a(t)}(\crnp)&= \Xi _+^{-a}e^+\overline
H_q^{t-a}(\rnp),\text{ for }t-a>-\tfrac1{q'},\\ 
C_*^{a(t)}(\crnp)&= \Xi _+^{-a}e^+\overline
C_*^{t-a}(\rnp),\text{ for }t-a>-1.
\endaligned
$$
 The $a$-transmission spaces are defined over $\Omega $ by localization.
When $\Omega $ is a bounded $C^{1+\tau }$-domain ($\tau
>0$), $u\in H_q^{a(t)}(\comega)$ is defined for $t<1+\tau $ to mean
that (1)--(2) hold:
\roster
\item $u$ is in $H^t_q$ on
compact subsets of $\Omega $.
\item Every $x_0\in\partial\Omega $ has an open neighborhood $U$ and a $C^{1+\tau }$-diffeomorphism in $\rn$ mapping $U'$
to $U$ such that $U'\cap \crnp$ is mapped to $U\cap \comega$, and $u$
is pulled back to a function $u'$  in
$H_q^{a(t)}(\crnp)$ locally (i.e., $\varphi u'\in H^{a(t)}(\crnp)$ when $\varphi
\in C_0^\infty (U')$).
\endroster
There is a similar definition of $C_*^{a(t)}(\comega)$.

A structural analysis as in Theorem 1.2 is valid also for  these
$a$-transmission spaces. Generally,  $\dot H^t_q(\comega)\subset H_q^{a(t)}(\comega) \subset \dot
H^a_q(\comega)$ for $t\ge a$, and there holds:
$$ 
H_q^{a(t)}(\comega)\cases =\dot H^t_q(\comega)\text{ when
  }-\tfrac1{q'}<t-a<\tfrac1q,\\
  \subset \dot H_q^{t(-\varepsilon )}
  (\overline\Omega )+d^ae^+\overline H_q^{t-a}(\Omega ),\text{ when } t-a
  >\tfrac1{q},
  \endcases \tag1.18
  $$
  where  $d(x)=\operatorname{dist}(x,\partial\Omega )$ near $\partial\Omega $
(extended positively to $\Omega $), and $(-\varepsilon )$ is active if $t-a-\frac1q$ is integer. The $d^a$-contribution can be described more exactly for specific values of
$t$. (Cf.\ \cite{G15}, \cite{G19}, \cite{AG23}.) Moreover, $\dot
H_q^t(\comega)\subset d^tL_q(\Omega )$ for $t\ge 0$, $\Omega $ smooth
(as kindly told us by Triebel, cf.\ e.g.\
 \cite{T12, Prop.\ 5.7}).   

There are similar statements for the $C^t_*$-scale; moreover,
  $C_*^{a(s+2a)}(\overline\Omega )\subset d^aC^{s+a}(\comega)$ when
  $s+2a,s+a\notin\Bbb N$, $s+a>0$, cf.\ \cite{G23}.

Note that in all the mentioned spaces with $t\ge a$, resp.\ $s>-a$,
the functions have a factor $d^a$ at the boundary.
\medskip

Let $P=(-\Delta )^a$ (or a suitable pseudodifferential
generalization explained further below),
and let $\Omega \subset \rn$ be a bounded
open set with some regularity. Recall that the homogeneous Dirichlet
problem is:
$$
r^+Pu= f \text{ in }\Omega ,\quad \operatorname{supp}u\subset
\overline\Omega .\tag1.3
$$
 
We know from the variational theory
that the Dirichlet realization $P_D$ in $L_2(\Omega )$ is bijective from $D(P_D)$ to
$L_2(\Omega )$, and ask now what can be said about $u$ when $f$ has
some regularity.  
Modern results:

\smallskip
$\bullet$       Ros-Oton and Serra showed in [RS14] by
potential-theoretic methods, when $\Omega $ is $C^{1,1}$:
$$
f\in L_\infty (\Omega )\implies u\in d^a\overline  C^t(\Omega ),\text{ for
  small $t>0$.}\tag1.19
$$
  The result
was extended later to $t$ up to $ a$. It was lifted to
higher-order H\"older spaces by Abatangelo and Ros-Oton in \cite{AR20}.

$\bullet$ The present author showed in \cite{G14},\cite{G15} by
pseudodifferential methods, when $\Omega $ is $C^\infty $, $1<q<\infty $:
$$
\align
f\in \overline H^s_q(\Omega )&\iff u\in
H_q^{a(s+2a)}(\overline\Omega),\text{ when }s >-a-1/q',\tag1.20 \\
  f\in \overline C_*^s(\Omega )&\iff u\in  C_*^{a(s+2a)}(\overline\Omega ),\text{ when }s >-a-1,\tag1.21\\
f\in C^\infty (\overline\Omega )&\iff u\in d^aC^\infty
                                  (\overline\Omega ).
                                   \tag1.22
                                  \endalign
$$

This theory initiated in an unpublished (and on some points sketchy)
lecture note of H\"ormander [H66]
(with $q=2$); (1.22) was obtained there.

  (1.20) is extended 
  to $C^{1+\tau }$-domains ($\tau >2a$) in a joint work with  Abels \cite{AG23},
  then valid for $0\le s<\tau -2a$.   The part $\implies$ in (1.21) is also obtained
  there with $s+2a$ replaced by $s+2a-\varepsilon $.

 Note the  {\it  sharpness} in (1.20)--(1.22);
they exhibit {\it the exact solution space} for (1.3). As pointed out above, the
functions there all have a factor $d^a$ near the boundary.

\example{Remark 1.3} An advantage of viewing $P$ as an elliptic pseudodifferential operator is that
we get {\it interior regularity} for free: When $f$ is locally in
$H^s_q$ (or $C^s_*$) in $\Omega $, then any solution of (1.3) is
locally in $H^{s+2a}_q$ (resp.\ $C^{s+2a}_*$) in $\Omega $. This has been
known since the advent of $\psi $do methods in the 1960's.
\endexample

 Now let us list the  hypotheses on general $\psi $do's
 $P=\operatorname{Op}(p(x,\xi ))$, under which our results hold.

\proclaim{Assumption 1.4} $P=\operatorname{Op}(p(x,\xi ))$ satisfies:

$1^\circ$ $p$ is {\bf classical} of order $2a>0$, i.e.,  $p\sim \sum_{j\in
  {\Bbb N}_0}p_j$ with $p_j(x,t\xi )=t^{2a-j}p_j(x,\xi )$ for $|\xi |\ge 1$.
The sign $\sim$ means that for all $J$, $\partial_x^\beta \partial_\xi ^\alpha [p-\sum_{j<J}p_j]$ is
$O(\ang\xi ^{2a-J-|\alpha |})$, for all multi-indices $\alpha ,\beta $.

$2^\circ$ $p$ is {\bf strongly elliptic:} $\operatorname{Re}p_0(x,\xi
)\ge c|\xi |^{2a}$ for $|\xi |\ge 1$, with $c>0$.
 
$3^\circ$ $p$ is {\bf even:}  $p_j(x,-\xi )=(-1)^jp_j(x,\xi )$, all
$j$,  $|\xi |\ge 1$.
\endproclaim

Assumption 1.4 is satisfied e.g.\ by $L^a$ when $L$ is a 2' order
strongly
elliptic differential  operator, and the $a$'th power is constructed as
in Seeley \cite{S67}, but also cases not stemming from differential
operators are included.

For a given smooth $\Omega $, it suffices for the results (1.20)--(1.22)
that $3^\circ$ holds for $p$ and derivatives $\partial_x^\beta \partial_\xi ^\alpha p$ at the points
$x\in\partial\Omega $, with $\xi $ just taken equal to the interior normal $\nu (x)$;
this is the so-called {\bf $a$-transmission condition} introduced by H\"ormander [H66],
[H85], also explained in \cite{G15}.

In \cite{AG23}, the hypotheses were generalized to allow symbols that are
only $C^\tau $ with respect to $x$, coupled with domains $\Omega $
that are only $C^{1+\tau }$; in this case (1.20) (and part of (1.21)) was  obtained
for $0\le s< \tau -2a$. 
\smallskip  
  
 Here are some words on the proof of (1.20), in the case where  $\Omega $ is
 $C^\infty $.  Roughly speaking, we perform two steps:

{\it  Step 1.} Reduce, by cut-downs and change-of-variables, to
situations where $\Omega $ is replaced by $\rnp$. Then $P$ is also modified.

  {\it   Step 2.} For the resulting $P$, let $Q=\Xi _-^{-a}P\Xi
  _+^{-a}$, so that 
$$
P=\Xi _-^{a}Q\Xi
  _+^{a}.\tag1.23
$$
  Here $Q$ is of order 0, and has some
  bijectivity properties (as a special case of an operator in the
  calculus of Boutet de Monvel \cite{B71}).
Namely. $r^+Qe^+$ is essentially
bijective
from $\overline H_q^t(\rnp)$ to itself for all $t\ge 0$.
Then
we find a solution operator $R= \Xi _+^{-a}e^+(r^+Qe^+)^{-1}(r^+\Xi _-^{-a}e^+)$,
 $$
 \overline H_q^s(\rnp)\overset {r^+\Xi _-^{-a}e^+}\to
 \longrightarrow \overline H_q^{s+a}(\rnp)
 \overset {(r^+Qe^+)^{-1}}\to
 \longrightarrow \overline H_q^{s+a}(\rnp)\overset{\Xi
   _+^{-a}e^+}\to
 \longrightarrow \Xi_+ ^{-a}e^+\overline H_q^{s+a}(\rnp),
$$
where the last space is the $a$-transmission space $ H_q^{a(s+2a)}(\crnp)$.

  The above explanation was simplified in particular on two points: 1) The $\Xi _\pm^t$
  should actually be
  replaced by a refined family $\Lambda _\pm^t$ with better
  pseudodifferential properties. 2) In some of the calculations, there is an error term of
  order $-\infty $ that has to be dealt with (a common feature of
  pseudodifferential calculations).

    Our proof for H\"older-Zygmund spaces follows the same lines,
    using that the pseudodifferential theory extends to such spaces.
    It also works for a wealth of other Besov- and Triebel-Lizorkin
    spaces, cf.\ \cite{G14}.

In the case of domains with finite smoothness, there was a need to expand the
(complicated) tools that exist for $\psi $do's with
nonsmooth $x$-dependence, cf.\ \cite{AG23}.

For Lipschitz domains (where the boundary is only $C^{0,1}$), there
are results about regularity and numerical methods e.g.\ by Acosta, Borthagaray and
Nochetto [AB17], [BN23], in basic spaces of Sobolev and  Besov types. There also exist studies where $f$ is given in spaces with powers of
$d$ as weights.

\pagebreak
 
\head 2. Further developments \endhead

\subhead 2.1. Evolution problems and resolvents \endsubhead

First we give a quick review of consequences of the analysis of $P_D$
for evolution problems (heat equations)
with homogeneous Dirichlet condition.
The basic problem is;
$$
\align
\partial_tu+r^+Pu&=f\text{ on }\Omega \times I ,\\
u&=0\text{ on }(\Bbb R^n\setminus\Omega )\times I, \tag2.1\\
u|_{t=0}&=0;
\endalign
$$ 
where $u$ and $f$ depend on $(x,t)$. Here  $I=\,]0,T[\,$ and $\Omega $
 is bounded, open and $C^{1+\tau }$ for suitable $\tau >0$; for simplicity we take zero initial data.
 
By Laplace transformation, the evolution problem is closely connected with the
stationary problem 
for $P-\lambda I$, where $\lambda \in \Bbb C$.

There is an easy
result in the $L_2$-framework:
Here $P_D$ is positive selfadjoint when $P=(-\Delta )^a$, and for more
general $P$ satisfying Assumption 1.4, $P_D$ is lower
semibounded with its discrete spectrum and numerical range contained in a sectorial region $$
M=\{\lambda \in\Bbb C\mid \operatorname{Re}\lambda +\beta \ge c_1>0,
|\operatorname{Im}\lambda |\le c_2(\operatorname{Re}\lambda +\beta )\}.$$
In particular, $\C\setminus M$ is in the resolvent set, and there is
a resolvent estimate 
$$
\|(P_D-\lambda )^{-1}\|_{\Cal L(L_2(\Omega ))}\le c_3\ang\lambda ^{-1}
\text{ for } \operatorname{Re}\lambda \le -\beta .\tag2.2
$$
Then standard old techniques show existence and uniqueness of a
 solution of (2.1) for $f\in L_2(\Omega \times I)$, and
$$
 f\in L_2(\Omega \times I) \iff u\in L_2(I;D(P_D))\cap \overline H^1( I;L_2(\Omega
))\text{ with }u(x,0)=0. \tag2.3
$$
Thanks to the analysis of $P_D$, we can in the right-hand side replace
$D(P_D)$ by $H^{a(2a)}(\overline\Omega )$, giving
a precise result. It is interesting that it only depends on $a$, not
on the value of the symbol $p$. (More details in \cite{G18a,b} for $\tau
=\infty $, \cite{G23} for $\tau >2a$.)
 
Now one can ask what happens if $f$ is in other spaces?

In the $L_2$-setting there is a
 functional analytic result from  Lions and Magenes' book [LM68] that
 can be applied to lift
(2.3) a small step in $x$ and a large step in $t$ [G18a,b],[G23]:
\smallskip

$\bullet$ For $k\in \Bbb N$,  $
r =\min \{2a, a+\frac12-\varepsilon \}$,
$$
\aligned
f\in L_2(I;\overline H^{r}(\Omega
))&\cap \dot H^k(\overline I;L_2(\Omega ))\implies \\
u &\in L_2(I; H^{a(2a+r)}(\overline\Omega ))\cap \overline H^{k+1}(
I;L_2(\Omega )).
\endaligned\tag2.4
$$

In $L_q$-spaces other techniques are needed. Here we have shown in
[G18a,b],[G23]:

$\bullet$ When $P$ satisfies Assumption 1.4 and is  $x$-independent and
symmetric, then for $1<q<\infty $,
$$
 f\in L_q(\Omega \times I) \iff u\in L_q(I;H_q^{a(2a)}(\overline\Omega ))\cap \overline H_q^1( I;L_q(\Omega
)) )\text{ with }u(x,0)=0.\tag2.5
$$
This is based on the fact that the $L_q$-Dirichlet realization
$P_{D,q}$
(whose domain satisfies $D(P_{D,q})=H_q^{a(2a)}(\overline\Omega )$)
is defined from a Dirichlet form in the sense of Fukushima,
Oshima and
Takeda \cite{FOT94} (also called sub-Markovian), allowing application of a result of
Lamberton \cite{L87}. This also implies an estimate like (2.2) with $L_2$
replaced by $L_q$. The time-regularity can then lifted by
use of general techniques of Amann \cite{A97}, and there are results for
other regularity classes with respect to $x$.

This type of solvability result is often called {\it maximal
$L_q$-regularity}, cf.\ e.g.\ Denk and Seiler \cite{DS15}. We expect that perturbation methods
would allow $x$-dependent symbols to some extent; there is work in
progress investigating this.

\smallskip
  
$\bullet$ In anisotropic
H\"older spaces 
$\overline C^{s,r }(\Omega \times I)=L_\infty
(I;\overline C^{s}(\Omega ))\cap L_\infty (\Omega ;\overline C^{r }( I
))$, Ros-Oton with coauthors Fernandez-Real and Vivas
[FR17], [RV18] have shown for $x$-independent symmetric operators, that
the regularity can be lifted as follows:
$$ 
f\in \overline C^{\gamma  , \gamma  /2a}(\Omega \times I)
\implies \partial_tu\in \overline C^{ \gamma  , \gamma  /2a}(\Omega
\times I'),\; u/d^a\in \overline C^{a+ \gamma  ,(a+ \gamma  )/2a}(\Omega
\times I'),
$$ 
when $  \overline {I'}\subset I$; here $\Omega $ is assumed
$C^{2+\gamma }$, and  $  0<\gamma <a $ with $
a+ \gamma  \notin{\Bbb N}$.
\smallskip

There have also been studies of evolution problems in numerical
analysis, e.g.\ by  Acosta, Bersetche and  Borthagaray [ABB19] in
$L_2$-Sobolev spaces over  Lipschitz domains. 
There is a very recent posting on results
in  $L_q$-Sobolev spaces weighted by powers of the distance $d(x)$ and
other functions, by Choi, Kim and Ryu [CKR23].

As another aspect, we mention that there is  an analysis  (in $C^\infty $-domains) [G19] showing that
the regularity of $u$ cannot be
lifted all the way to $C^\infty(\comega\times \bar I) $
 or $d^aC^\infty(\comega\times \bar I) $ when $f\in
 C^\infty(\comega\times \bar I)$. This is  in contrast
with heat problems for the local operator $\Delta $.

  \subhead 2.2. Motivation for local nonhomogeneous boundary  conditions \endsubhead

Now we turn to nonhomogeneous Dirichlet conditions \cite{G15}, which will be
  explained in detail.
  
  As a nonhomogeneous Dirichlet problem, much of the literature
  considers the problem
$$
r^+Pu=f\text{ in }\Omega , \quad u=g \text{ on }\rn\setminus \Omega
,\tag2.6
$$
where the difference from (1.3) is that $u$ may take a nonzero value $g$
outside of $\Omega $.
 
There is an easy reduction of this problem to
the homogeneous case, namely: Let $G$ be a function extending $g$ to
$\rn$, then the problem (2.6) can be turned into the homogeneous
problem
$$
r^+Pu'=f'\text{ in }\Omega , \quad u'=0 \text{ on }\rn\setminus \Omega
,\tag2.7
$$
where $u'=u-G$, $f'=f-r^+PG$. The discussion of regularity of
solutions then involves how the extension from $g$ to $G$ is performed
and how it influences $r^+PG$.

We shall here discuss another Dirichlet condition that involves a
boundary value on $\partial\Omega $ and is {\it local}.
 For the motivation, consider $C^\infty $-results.
Define for any $\mu >-1$:
$$
\Cal E_\mu (\comega)=e^+d^\mu C^\infty (\comega).
$$
(As usual, $e^+ $ means extension by zero.) Here
$\Cal E_0(\comega)\simeq C^\infty (\comega)$.

With this notation, the regularity result (1.22) for $(-\Delta )^a$
and for the
generalizations $P$ satifying Assumption 1.4 states that
$$
f\in C^\infty (\overline\Omega )\iff u\in \Cal E_a (\overline\Omega
                                  ). \tag2.8
$$
Moreover, one can show the forward mapping property for all integers $k\ge -1$ [G15]
$$
r^+P\colon \Cal E_{a+k}(\comega)\to C^\infty(\comega).
$$
 
There are Taylor expansions at the boundary, in local coordinates where $\Omega
$ is replaced by $\rnp=\{x=(x',x_n)\mid x_n>0\}$ so that $d(x)=x_n$:
$$
\aligned
        &\text{In }\Cal E_0: u(x)\sim
v_0(x')+v_1(x')x_n+v_2(x')x_n^2+\dots,\quad  \text{ when }x_n>0.\\
&\text{In }\Cal E_1: u(x)\sim v_0(x')x_n+ v_1(x')x_n^2+v_2(x')x_n^3+\dots.\\
&\text{In }\Cal E_a: u(x)\sim v_0(x')x_n^a+v_1(x')x_n^{a+1}+
v_2(x')x_n^{a+2}+\dots.\\
&\text{In }\Cal E_{a-1}: u(x)\sim v_{0}(x')x_n^{a-1}+v_1(x')x_n^{a}+ v_2(x')x_n^{a+1}+\dots.
\endaligned
$$ 
Recall the notation
 $u|_{\partial\Omega }=\gamma _0u$. Note that the expansions of functions in $\Cal
 E_{a-1}$ only differ from those in $\Cal E_a$ by having a term
 $v_{0}(x')x_n^{a-1}$; i.e., $\gamma _0(u/x_n^{a-1})$ can be
 nontrivial. This leads to the important
 observation:
  $$ 
\Cal E_a\text{ is the subset of $\Cal E_{a-1}$ where }\gamma
_0(u/d^{a-1})=0.\tag2.9
$$
(It also holds when $a \ge 1$.)

  Let $f\in C^\infty
 (\overline\Omega )$, $\varphi \in C^\infty (\partial\Omega
 )$, for a bounded $C^\infty $-domain $\Omega $, and let us compare boundary value problems for $\Delta $ and $(-\Delta )^a$:

  \smallskip  
{\it Old fact:}   The {\bf nonhomogeneous Dirichlet problem for $\Delta $:}
$$
\aligned
\Delta u&=f \text{ on }\Omega ,\\
 \gamma _0u&=\varphi \text{ on
}\partial\Omega ,
\endaligned \tag 2.10
$$
is uniquely  solvable in $C^\infty (\overline\Omega )\simeq \Cal
E_0(\comega )$.  
\smallskip

 As a special case, the {\bf homogeneous Dirichlet problem for $\Delta $:} 
$$
\aligned
\Delta u&=f \text{ on }\Omega , \\
\gamma _0u&=0 \text{ on
}\partial\Omega ,
\endaligned \tag 2.11
$$
is uniquely  solvable in $\{ u\in C^\infty (\overline\Omega )\mid
\gamma _0u=0\}\simeq \Cal
E_1(\comega )$, cf.\ also (2.9).
\smallskip 

{\it Modern result:}  The {\bf homogeneous Dirichlet problem for $(-\Delta )^a$}
$$
\aligned
(-\Delta )^au&=f \text{ on }\Omega , \\
\operatorname{supp}u&\subset \overline\Omega ,
\endaligned 
$$
is uniquely  solvable in $ \Cal
E_a(\comega )$ (as already stated in (1.22) and (2.8)). 
 Here $\Cal E_a(\overline\Omega )$ has a role like the one $\Cal
 E_1(\overline\Omega )$ has 
for $\Delta$. 
 
\smallskip

Now it is natural to define a {\bf nonhomogeneous Dirichlet problem for
  $(-\Delta )^a$} by  going
out to the larger space $\Cal E_{a-1}(\overline\Omega )$.
The problem
$$
\aligned
(-\Delta )^au&=f \text{ on }\Omega , \\
\gamma _0(u/d^{a-1})&=\varphi  \text{ on
}\partial\Omega ,\\
\operatorname{supp}u&\subset \overline\Omega ,
\endaligned \tag2.12
$$
is uniquely  solvable in $\Cal
E_{a-1}(\comega )$. ({\it Proof:} subtract a function $w\in \Cal E_{a-1}$
with $\gamma _0(w/d^{a-1})=\varphi $, then $v=u-w$  solves a
homogeneous Dirichlet problem, cf.\ (2.8), (2.9).)
\smallskip

This is surprisingly simple! It can be generalized to solvability
statements in Sobolev spaces after some more work; see later.

The interest of the nonhomogeneous Dirichlet problem (2.12) was also pointed
out by Abatangelo [A15], from a very different viewpoint: He started
with a Green's function $G_\Omega (x,y)$ for the homogeneous Dirichlet
problem for $(-\Delta )^a$, and developed integral representation formulas imitating the
formulas known for $\Delta $, arriving at a strange boundary
operator $u\mapsto Eu$, that he showed was proportional to $\gamma
_0(u/d^{a-1})$ in the case where $\Omega $ is a ball.

$E$ is defined
by an integral formula; a proof that $Eu= c_0 \gamma _0(u/d^{a-1})$
for more general $\Omega $ is given 
in [G23] (the constant $c_0$ equals $\Gamma (a)\Gamma (a+1)$).
 This boundary operator also enters in other studies, e.g.\ by Chan, Gomez-Castro and
 Vazquez [CGV21], and by Fernandez-Real and Ros-Oton \cite{FR20}.

The solutions in $\Cal E_{a-1}$ are generally unbounded
on $\Omega $, since $u$ behaves like the unbounded function $d^{a-1}$
near $\partial\Omega $ (when $\varphi \ne 0$). They are therefore
often called ``blow-up solutions''. They are in $L_q(\Omega )$ for
$q<(1-a)^{-1}$.

There is also a local {\it Neumann condition} $\gamma _1(u/d^{a-1})=\psi $,
which has a good solvability theory \cite{G14},\cite{G18}; here $\gamma _1v=\gamma
_0\partial_\nu v$, the normal derivative.

\smallskip

For solvability results in general Sobolev spaces, the
role of $\Cal E_{a-1}$ will for $\rnp$ be taken over by the
$(a-1)$-transmission spaces defined  by
$$
H^{(a-1)(t)}(\crnp)=\Xi _+^{-a+1}e^+\overline H^{t-a+1}(\rnp),
$$
and $L_q$-variants with $q\ne 2$.
The model problem (as in Example 2 above) is now: 
$$
\aligned
r^+(1-\Delta )^au&=f\text{ in }\rnp ,
\\
\gamma _0(u/x_n^{a-1})&=\varphi \text{ on }\Bbb R^{n-1} ,\\
\operatorname{supp}u&\subset \crnp.
\endaligned \tag2.13
$$
 
First we observe that there is a result on boundary values like in Theorem
1.2 but with $a$ replaced by $a-1$:

\proclaim{Theorem 2.1}
The mapping $\gamma _0^{a-1}\colon u\mapsto \gamma
  _0(u/x_n^{a-1})$  from $x_n^{a-1}\Cal S(\crnp)$ to $\Cal S(\R^{n-1})$ extends to a continuous surjective
  mapping 
(when $t>a-\frac12$),
$$
\gamma _0^{a-1}\colon H^{(a-1) (t)}(\crnp )\to H^{t-a+\frac12}(\Bbb
  R^{n-1} ).
  $$
Here  $H^{a(t)}(\crnp )$ is a closed subspace of $H^{(a-1) (t)}(\crnp )$,
equal to the set where $\gamma _0^{a-1}u=0$.
\endproclaim 

 The last line comes from (2.9).

Then we solve (2.13) by subtracting from $u$ a term $w$ with $\gamma
_0^{a-1}w=\varphi $, reducing to the homogeneous Dirichlet problem.
As a result (note that $s+2a$ plays the role of $t$):

  \proclaim{Theorem 2.2}
 The nonhomogeneous Dirichlet problem {\rm (2.13)}
with given $f\in \overline H^s(\rnp )$, $\varphi \in
H^{s+a+\frac12}(\Bbb R^{n-1})$, $s \ge 0$,
is uniquely solvable with a
solution $u\in H^{(a-1)(s+2a)}(\crnp)$.
\endproclaim

\subhead 2.3 Nonhomogeneous Dirichlet conditions over curved
domains \endsubhead

For curved domains $\Omega $, the $(a-1)$-transmission spaces are defined by use of local coordinates.
For the $H^s_q$-scales with $q\ne 2$, the correct spaces over the
boundary are Besov spaces $B^t_{q}$ (also denoted  $B^t_{q,q}$). Here
the trace map $\gamma _0^{a-1}u=\gamma _0(u/d^{a-1})$ satisfies that
 
$$\gamma _0^{a-1} \colon  H_q^{(a-1) (t)}(\comega )\to
B_q^{t-a+\frac1{q'}}(\partial \Omega )$$
is continuous and surjective for
$t>a-\frac1{q'}$, with kernel  $H_q^{a(t)}(\comega )$. One finds:

 \proclaim{Theorem 2.3} 
There is unique solvability of the nonhomogeneous Dirichlet problem
$$
\aligned
Pu&=f\text{ in }\Omega ,
\\ \gamma _0^{a-1}u&=\varphi \text{ on }\partial\Omega ,\\
\operatorname{supp}u&\subset \comega,
\endaligned \tag2.14
$$
for given $f\in \overline H_q^s(\Omega )$, $\varphi \in
B_q^{s+a+1/q'}(\partial\Omega )$, $s\ge 0$, with 
solution $u\in H_q^{(a-1)(s+2a)}(\comega)$. 
\endproclaim

This is shown is [G15] for bounded smooth $\Omega $, under Assumption 1.4.
 (More precisely, if $P\ne (-\Delta )^a$, $0$ can be an eigenvalue of
 the homogeneous Dirichlet problem, and in that case, there is only a
 Fredholm  solvability.) In \cite{G23} the result is generalized to  $C^{1+\tau
 }$-domains $\Omega $ and $\psi $do's $P$ with $C^\tau $ $x$-dependence,
 when $0\le s<\tau -2a-1$.

\smallskip 

 These stationary results can be followed up with results for evolution problems (for
  $\partial_t+P$) and
  resolvent problems (for $P-\lambda $, $\lambda \in\Bbb C$):

For the study of (2.14) with $P$ replaced by $P-\lambda $, we need $u$
to be at least in $L_q(\Omega )$. The domain space
$H_q^{(a-1)(s+2a)}(\comega)$ ($s\ge 0$) is not always there. In fact,
already for $s=0$ (recall $1<q<\infty $),
$$
 H_q^{(a-1)(2a)}(\comega)\subset L_q(\Omega ) \text{ if and only if }q<(1-a)^{-1}.\tag2.15
 $$
  (For $q=2$, this holds when $a>\frac12$.) 

The evolution problem is:
$$ 
\aligned
Pu+\partial_tu&=f\text{ on }\Omega \times I ,\\
u&=0\text{ on }(\Bbb R^n\setminus\Omega )\times I, \\
\gamma _0^{a-1}u&=\psi \text{ on }\partial\Omega \times I,\\
u|_{t=0}&=0.
\endaligned\tag2.16
$$

Here we can show \cite{G23}:
\proclaim{Theorem 2.4}
 Let $q<(1-a)^{-1}$. If $q\ne 2$, let $P$ be $x$-independent  symmetric.   
  For $f(x,t)$ given in $L_q(\Omega \times I)$, and $\psi (x,t)$ given
  in $ L_q(I;
B_q^{a+1/{q'}}(\partial\Omega ))\cap \ol H_q^1(I; B_q^{\varepsilon
}(\partial\Omega ))$ 
 with $\psi (x,0)=0$ (some $\varepsilon >0$), there is a
  unique solution $u(x,t)$ of {\rm (2.16)} satisfying
  $$
u\in L_q(I;H_q^{(a-1)(2a)}(\comega ))\cap \overline H_q^1(I; L_q(\Omega )).
$$
\endproclaim

It is shown by reduction to a problem with $\psi =0$, where
(2.3)--(2.5) can be applied.

Solvability of resolvent problems is obtained in the following
theorem \cite{G23, Th.\ 5.4}:

\proclaim{Theorem 2.5}  Let $q<(1-a)^{-1}$. If $q\ne 2$, let $P$ be
$x$-independent  symmetric. Denote by $\Sigma $ the spectrum of $P_D$
(it is discrete). Consider for $\lambda \in\C$ the problem
 $$
 \aligned
     Pu-\lambda u& = f \text { in }\Omega ,\\
     u&=0 \text { in }\rn\setminus \Omega ,\\
     \gamma _0^{a-1}u&=\varphi \text{ on }\partial\Omega ,
     \endaligned \tag2.17
$$
with $f$
given in $L_q(\Omega )$, $\varphi $ given in $
B_q^{a+1/q'}(\partial\Omega )$, and the solution being sought in 
$H_q^{(a-1)(2a)}(\comega)$.

If $\lambda
\notin\Sigma $, it is uniquely solvable.

If $\lambda \in\Sigma $, it is Fredholm solvable,  with the same
dimension of the kernel and cokernel of the mapping $u\to \{f,\varphi
\}$.
\endproclaim 

There are related resolvent studies  by Chan, Gomez-Castro and Vazques [CGV21]
in weighted $L_1$-spaces, generally larger than the spaces we consider
in \cite{G23}.  For $f=0$, \cite{CGV21} regards (2.17) as an
"eigenvalue problem", and presents it as a mysterious fact that the
solutions ("eigenfunctions") generally blow up at the boundary. We
find this natural, since the functions
in the precise domain $H_q^{(a-1)(2a)}(\comega)$ have a factor $d^{a-1}$ at the boundary
as soon as $\gamma _0(u/d^{a-1}) $ is nontrivial.

\subhead 2.4. Integration by parts,  Green's formula \endsubhead

Another topic that we shall touch upon very briefly is the question of  integration by parts formulas for the
 fractional Laplacian and its generalizations.
Ros-Oton and Serra [RS14a] started the analysis by showing a Pohozaev formula
for solutions of the homogeneous Dirichlet problem, important for
uniqueness questions in nonlinear variants.  
Their basic result is, in an equivalent version: 

\proclaim{Theorem 2.6}
 Let $\Omega $ be bounded and $C^{1,1}$. Let $u$
  and $v$ be solutions of the homogeneous Dirichlet problem {\rm
 (1.3)} for
  $(-\Delta )^a$ with real
  right-hand side in $L_\infty (\Omega )$, so they are in
  $d^aC^t(\comega)$ (small $t$) by {\rm (1.19)}. Then for each $j$,
  $$
\int_\Omega ((-\Delta )^au \,\partial_j v+
\partial_ju\,(-\Delta )^a v)\, dx =\Gamma (a+1)^2\int_{\partial\Omega }\nu
_j\gamma _0(\tfrac u{d^a})\,\gamma _0(\tfrac v{d^{a}})\, d\sigma ,\tag 2.18
$$
where $\nu =(\nu _1,\dots,\nu _n)$ is the interior normal.
\endproclaim 

Their proof is based on a fine analysis of the factorization $(-\Delta
)^a=(-\Delta )^{a/2}(-\Delta )^{a/2}$ applied to real functions. In \cite{G16}, we
worked out a proof  of (2.18) based on Fourier analysis and factorizations
developed from
(1.23), applicable to operators satisfying Assumption 1.4
and smooth domains.

Moreover, we have shown integration formulas also for 
solutions of
nonhomogeneous boundary problems.
Let us go directly to  the Green's formula \cite{G18}, \cite{G20}:

\proclaim{Theorem 2.7}
 Let $\Omega $ be bounded smooth. For $u,v\in
  H^{(a-1)(s)}(\comega)$ there holds when $s>a+\frac12$:
$$
\int_\Omega \bigl((-\Delta )^au \,\bar v-
u\,(-\Delta )^a\bar v\bigr)\, dx =c_0\int_{\partial\Omega }\bigl(\gamma
_1(\tfrac u{d^{a-1}})\,\gamma _0(\tfrac{\bar v}{d^{a-1}})
-\gamma _0(\tfrac u{d^{a-1}})\,\gamma _1(\tfrac{\bar v}{d^{a-1}})\bigr)\, d\sigma ,
\tag2.19$$
$c_0= \Gamma (a)\Gamma
(a+1)$.  
\endproclaim

Note that both the Dirichlet trace   $\gamma
_0(\tfrac u{d^{a-1}})$ and the Neumann trace   $\gamma
_0(\partial_\nu (\tfrac u{d^{a-1}}))$ enter in (2.19). When $\gamma
_0(\tfrac u{d^{a-1}})=0$, the Neumann trace equals the value $\gamma _0(\tfrac u{d^a})$ entering in (2.18).

For general $P$ satisfying Assumption 1.4, there is a similar formula with an extra term
$\int_{\partial\Omega }B\gamma _0^{a-1}u\,\gamma _0^{a-1}\bar v\,dx$,
where $B$ is a $\psi $do on $\partial\Omega $ of
order 1.  
\medskip

  We end this survey by some remarks on what more can be done, or
  needs doing, in the present context. Here are a few suggestions:
  \roster  
\item More on  evolution problems  in $L_p$-Sobolev spaces, also for
$x$-dependent opera\-tors $P$.
\item Development from \cite{G14} of consequences in $L_1$-spaces and
in general $F^s_{p,q}$- and
$B^s_{p,q}$-spaces. 
\item Extension of more results known for smooth domains (e.g.\
integration formulas), to nonsmooth domains. 
\item Applications to problems with  nonlinearity.
\item Treatment of operators without the reflection symmetries of $(-\Delta )^a$. 
\endroster

Ad (5):  Ros-Oton and colleagues have initiated  studies of boundary
  value problems for operators
  that do not have the {\it evenness} property of  $(-\Delta )^a$ and the
  operators $P$ we have listed. For
  example $(-\Delta )^\frac12+b\cdot \nabla$, $b\in\rn$, with an even
 part $(-\Delta )^\frac12$ and an odd part $b\cdot \nabla$. They get
  results by real integral operator methods (from potential theory and function
  theory); for a comprehensive treatment see  Dipierro, Ros-Oton, Serra and Valdinoci
  [DRSV22].

 By Fourier methods we can treat completely general strongly elliptic operators $L=\operatorname{Op}(\ell(\xi ))$,
 where $\ell(\xi )$ is homogeneous of order $2a$ and just satisfies
 $\operatorname{Re}\ell(\xi )\ge c|\xi |^{2a}$ with $c>0$, showing how
 a $\mu $-transmission space comes in (with a possibly complex $\mu $), and obtaining an integration
 by parts formula; but so far
 only in the model case of $\rnp$ [G22]. It might be worth trying to apply the localization techniques of [DRSV22] to
 extend the results for $L$ to curved domains.

\enddocument